\documentclass[a4,12pto] {article}
\author{Mar\'{\i}a de la Paz Tirado Hern\'andez}
\usepackage{amssymb}
\usepackage{amsfonts}
\usepackage{latexsym}

\usepackage{amsmath}
\usepackage[latin1]{inputenc}

\usepackage[all]{xy}

\usepackage{theorem}
\usepackage{graphicx}
\usepackage{pb-diagram}
\usepackage{verbatim}

\usepackage{anysize}
\usepackage{cancel}
\usepackage{color}
\usepackage[colorlinks=true,linkcolor=blue,urlcolor=black]{hyperref}

\newtheorem{teo}{Theorem}[section]
\newtheorem{prop}[teo]{Proposition}
\newtheorem{cor}[teo]{Corollary}
\newtheorem{lem}[teo]{Lemma}
\newtheorem{Def}[teo]{Definition}
\newtheorem{nota}[teo]{Remark}

\newtheorem{ej}[teo]{Examples}

\newcommand{\h}{\langle h \rangle}

\DeclareMathOperator{\Ann}{\rm Ann}

\DeclareMathOperator{\val}{\rm val}

\DeclareMathOperator{\IDer}{\rm IDer}

\DeclareMathOperator{\Der}{\rm Der}

\DeclareMathOperator{\HS}{\rm HS}

\DeclareMathOperator{\Id}{\rm Id}

\usepackage{fancyhdr}
\marginsize{2cm}{2cm}{2.5cm}{2.5cm}

\title{ Integrable derivations in the sense of Hasse-Schmidt for some binomial plane curves}
\author{Mar\'{\i}a de la Paz Tirado Hern\'andez\thanks{Partially supported by MTM2016-75027, P12-FQM-2696 and FEDER.} \thanks{Departamento de \'Algebra e Instituto de Matem\'aticas (IMUS), Universidad de Sevilla, Espa\~na.}}

\begin{document}
\date {}
\maketitle
\begin{abstract}
We describe the module of integrable derivations in the sense of Hasse-Schmidt of the quotient of the polinomial ring in two variables over an ideal generated by the equation $x^n-y^q$.


\end{abstract}

\begin{center} INTRODUCTION \end{center}

Let $k$ be a commutative ring and $A$ a commutative $k$-algebra. A Hasse-Schmidt derivation of $A$ over $k$  of length $m\in \mathbb N$ or $m=\infty$ is a sequence $D=(D_i)_{i\geq 0}^m$ such that: $$
\begin{array}{ccc}
D_0=\Id_A,&\displaystyle D_i(xy)=\sum_{a+b=n} D_a(x)D_b(y)
\end{array}
$$
for all $x,y\in A$. We denote by  $\HS_k(A;m)$ the set of
Hasse-Schmidt derivations of $A$ of length $m$. The component $D_i$ of a Hasse-Schmidt derivation is a differential operator of order $\leq i$, in particular $D_1$ is a $k$-derivation.

The Hasse-Schmidt derivations of length $m$, also called higher derivation of order $m$ (see \cite{Ma}), were introduced by H.Hasse and F.K. Schmidt (\cite{H-S}) and they have been used by several authors in different contexts (see \cite{Na1}, \cite{Se} or \cite{Tr}). An important notion related with Hasse-Schmidt derivations is integrability. Let $m\in \mathbb N$ or $m=\infty$, then we say that
$\delta\in \Der_k(A)$ is $m$-integrable if there exists $D\in
\HS_k(A;m)$ such that $\delta=D_1$. The set of all $m$-integrable $k$-derivations is an $A$-submodule of $\Der_k(A)$ for all $m$, which is denoted by $\IDer_k(A;m)$.

If $k$ has characteristic 0 or $A$ is $0$-smooth over $k$, then any $k$-derivation is $\infty$-integrable (\cite{Ma}), that means that $\Der_k(A)=\IDer_k(A;\infty)$. If we consider $k$ a ring of positive characteristic and $A$ any commutative $k$-algebra, the modules $\IDer_k(A;m)$ have better properties than $\Der_k(A)$ (see \cite{Mo}). So exploring these modules seems interesting to better understand singularities in positive characteristic.


The aim of this paper is to describe the modules of $m$-integrable derivations, for $m\geq 1$ and $m=\infty$, of the quotient of the polynomial ring in two variables over an ideal generated by an equation of type
$x^n-y^q$.

This paper is organized as follows:
In section 1 we recall the definition of Hasse-Schmidt derivations and give some known properties that will be useful in later sections. In section 2 we focus on the integrability of derivations in the sense of Hasse-Schmidt in quotients of polynomial rings in two variables over the ideal generated by the equation $x^n-y^q$. Namely, we calculate the module of integrable $k$-derivations when $k$ is a reduced ring of characteristic $p>0$ and  $n$ or $q$ are not multiple of $p$. In section 2.1, we assume that $k$ is a unique factorization domain and we see the relationship between integrable derivations of the quotient of a polynomial ring over $\langle f\rangle$ and over $\langle f^p\rangle$ where $f$ is a polynomial. Thanks to this relationship, we can describe the integrable derivations of $k[x,y]/\langle x^n-y^q\rangle$ when $n$ and $q$ are both multiples of $p$. In section 3, we calculate the module of integrable derivations in some examples taken from \cite{Gr}.

{\bf Acknowledgment.} The author thanks Professor Luis Narv\'aez Macarro for their careful reading of this paper with numerous useful comments.

\section{Hasse-Schmidt derivations}
Let $k$ be any commutative ring and $A$ a commutative $k$-algebra. In this section we will define Hasse-Schmidt derivations and we will give some of their properties, ending with the case where A is a polynomial ring. We denote $\overline{\mathbb
N}=\mathbb N \cup \{\infty\}$. For each integer $m\geq 1$, we will write $A[|\mu|]_m:=A[|\mu|]/\langle \mu^{m+1}\rangle$ and
$A[|\mu|]_\infty:=A[|\mu|]$.

\begin{Def}\label{DefHS}
A Hasse-Schmidt derivation (over $k$) of $A$ of length $m\geq 1$ (resp. of
length $\infty$) is a sequence $D:=(D_0,D_1,\ldots, D_m)$
(or resp. $D=(D_0,D_1,\ldots)$) of $k$-linear maps $D_i:A\rightarrow
A$, satisfying the conditions:
$$
\begin{array}{ccc}
D_0=\Id_A,&\displaystyle D_i(xy)=\sum_{a+b=n} D_a(x)D_b(y)
\end{array}
$$
for all $x,y\in A$ and for all $i$. We write $\HS_k(A;m)$ (resp.
$\HS_k(A)$) for the set of Hasse-Schmidt derivations (over $k$) of $A$ of
length $m$ (resp. $\infty$).
\end{Def}

\begin{nota}[\cite{Ma}; cf. \cite{Na2}]

\begin{enumerate}
\item Any Hasse-Schmidt derivation $D\in \HS_k(A;m)$ is determined by the $k$-algebra homomorphism
$$
\begin{array}{rccc}
\varphi_D:&A&\rightarrow& A[|\mu|]_m \\ &a&\mapsto& \displaystyle \sum_{i\geq 0}^m
D_i(a)\mu^i
\end{array}
$$
satisfying $\varphi_D(x)=x\mod \mu$. $\varphi_D$ can be uniquely
extended to a $k$-algebra automorphism
$\widetilde\varphi_D:A[|\mu|]_m\rightarrow A[|\mu|]_m$ with
$\widetilde\varphi_D(\mu)=\mu$. So, $\HS_k(A;m)$ has a canonical
group structure. Namely, $D\circ D'=D''\in \HS_k(A;m)$ with
$D''_n=\sum_{i+j=n} D_i\circ D_j'$ for $n\leq m$. Moreover, the
component $D_1$ is a $k$-derivation. So, the map $(\Id, D_1)\in
\HS_k(A;1)\mapsto D_1\in \Der_k(A)$ is a group isomorphism.

\item \label{Multporesc} For any $a\in A$ and any $D\in \HS_k(A;m)$, the sequence $a\bullet D=(a^iD_i)\in
\HS_k(A;m)$.

\item For any $1\leq n\leq m$ and any $D\in \HS_k(A;m)$, we define the truncation map by $\tau_{mn}(D)=(\Id,D_1,\ldots,D_n)\in \HS_k(A;n)$.
\end{enumerate}
\end{nota}

\begin{Def}\label{Log-Int} Let $D\in \HS_k(A;m)$ where $m\in \overline{\mathbb N}$ and $n\geq m$. Let $I$ be an ideal of $A$.
\begin{itemize}
\item We say that $D$ is $I$-logarithmic if $D_i(I)\subseteq I$ for all $i$. The set of $I$-logarithmic Hasse-Schmidt derivations is denoted by $\HS_k(\log I;m)$, $\HS_k(\log I):=\HS_k(\log I;\infty)$ and
$\Der_k(\log I):=\HS_k(\log I;1)$.
\item We say that $D$ is $n$-integrable if there exists $E\in \HS_k(A,n)$ such that $\tau_{nm}(E)=D$. Any such $E$
will be called a $n$-integral of $D$. If $D$ is $\infty$-integrable
we say that $D$ is integrable. If $m=1$, we write $\IDer_k(A;n)$ for
the set of $n$-integrable derivations and
$\IDer_k(A):=\IDer_k(A;\infty)$.
\item We say that $D$ is $I$-logarithmically $n$-integrable if there exists $E\in \HS_k(\log I;n)$ such that $E$ is a $n$-integral of
$D$. We put $\IDer_k(\log I;n)$ for the set of $I$-logarithmically $n$-integrable derivations  when $m=1$ and
$\IDer_k(\log I):=\IDer_k(\log I, \infty)$.
\end{itemize}
\end{Def}

\begin{nota}
$\IDer_k(A;n)$ is an $A$-submodule of $\Der_k(A)$ thanks to the
group structure of $\HS_k(A;n)$ and operation \ref{Multporesc}.
\end{nota}

\begin{Def}\label{leap}
$A$ has a leap on $s>1$ if the inclusion $\IDer_k(A;s-1)\supsetneq
\IDer_k(A;s)$ is proper.
\end{Def}

\begin{lem}\label{Dihp}
Let $k$ be a ring of characteristic $p>0$ and $h\in A$. Consider
$D\in \HS_k(A;m)$ with $m\in \overline{\mathbb N}$ and $\tau\geq 0$.
Then, for all $i\leq m$,
$$
D_i\left(h^{p^\tau}\right)=\left\{\begin{array}{lcl}
0&\mbox{ if }& p^\tau \not| i\\
D_{i/p^\tau}(h)^{p^\tau}& \mbox{ if }&p^\tau|i
\end{array}\right.
$$
\end{lem}

\noindent{\bf Proof.}

Let $\varphi:A\rightarrow A[|\mu|]_m$ be the $k$-algebra
homomorphism determined by $D$. Then,
$$
\sum_{i\geq 0}^m
D_i\left(h^{p^\tau}\right)\mu^i=\varphi\left(h^{p^\tau}\right)=
\varphi(h)^{p^\tau}=\sum_{j\geq 0}^m D_j(h)^{p^\tau}\mu^{jp^\tau}
\mod \left\langle \mu^{m+1}\right\rangle
$$
and we obtain the result by equating the coefficients in the  above equation.
\vspace{-0.5cm}\begin{flushright}$\square$\end{flushright}

\begin{lem}\label{ExpDi}
Consider $g\in A$ and $D\in \HS_k(A;m)$. Suppose that $D_j(g)\in
\langle g\rangle$ for all $0\leq j<m$. Then, for all $r\geq 1$,
$$
D_m(g^r)\in rg^{r-1}D_m(g)+\langle g^r\rangle.
$$
\end{lem}

\noindent{\bf Proof.}

We will prove that $D_j(g^r)\in \langle g^r\rangle$ for all $j<m$
and $r\geq 1$. We proceed by induction on $j$. For $j=0$ the result
is clear since $D_0=\Id$. Let us assume that $D_a(g^r)\in \langle
g^r\rangle$ for all $a<j$ and all $r$. We will show the result for
$j$ by induction on $r$. When $r=1$, it's obvious from the
hypothesis. Let us suppose that $D_j(g^{r-1})\in \langle
g^{r-1}\rangle$. From the definition of Hasse-Schmidt derivation,
$$
D_j\left(g^r\right)= D_j\left(g^{r-1}\right) g+\sum_{\substack{
a+b=j\\a,b\neq 0}} D_a\left(g^{r-1}\right) D_b(g)+g^{r-1}D_j(g)\in
\langle g^r\rangle.
$$
Now, we will prove the lemma by induction on $r\geq 1$. It is
obvious for $r=1$, let us suppose that $D_m(g^{r-1})\in
(r-1)g^{r-2}D_m(g)+\langle g^{r-1}\rangle$. From the definition of
Hasse-Schmidt derivation,
$$
D_m\left(g^r\right)=D_m\left(g^{r-1}\right)g+D_m(g)g^{r-1}+\sum_{\substack{
a+b=m\\a,b\neq 0}} D_a\left(g^{r-1}\right)D_b(g) \in
rg^{r-1}D_m(g)+\left\langle g^r\right\rangle
$$
and the lemma is proved.
\vspace{-0.7cm}\begin{flushright}$\square$\end{flushright}

\subsection{Polynomial ring and integrability}

Consider $R=k[x_1,\ldots,x_d]$ the polynomial ring over a commutative ring $k$. In this section, we recall, for the ease of the reader, some results related with the integrability of $k$-derivation in a polynomial ring.

\begin{teo}{\rm\cite[{\bf Th. 27.1}]{Ma}}\label{IDer=Der}
Let $R=k[x_1,\ldots,x_d]$ the polynomial ring over $k$, then $\IDer_k(R)=\Der_k(R)$.
\end{teo}

\begin{cor}\label{HSenpolise extiende}
Any Hasse-Schmidt derivation of $R$ over $k$ of length $m\geq 1$ is
integrable.
\end{cor}

\noindent{\bf Proof.} This is consequence of Theorem \ref{IDer=Der}
and Proposition 2.1.5 of \cite{Na2}.
\vspace{-0.7cm}\begin{flushright}$\square$\end{flushright}

\begin{cor}{\rm\cite[{\bf Corollary. 2.1.10}]{Na2}}\label{apli}
The map $\Pi:\IDer_k(\log I;m)\rightarrow \IDer_k(R/I;m)$ defined by
$\Pi(D)=\overline{D}$ where $\overline{D_i}(a+I)=D_i(a)+I$ is a
surjective group homomorphism.
\end{cor}

\begin{cor}\label{Saltos del cociente}
Let $I\subset R$ be an ideal and $A=R/I$. Then, $A$ has a leap on
$s\geq 1$ if and only if the inclusion $\IDer_k(\log I;s-1)\supsetneq
\IDer_k(\log I;s)$ is proper.
\end{cor}

\begin{prop}{\rm\cite[{\bf Prop. 2.2.4}]{Na2}} \label{2.2.4Na2}
Let $f\in R$, $I=\langle f\rangle$, and $J^0=\langle \partial_1(f),\ldots,\partial_d(f)\rangle$ the gradient ideal. If $\delta:R\rightarrow R$ is an $I$-logarithmic $k$-derivation with $\delta\in J^0\Der_k(R)$, then $\delta$ admits an $I$-logarithmic integral $D\in \HS_k(\log I)$ with $D_i(f)=0$ for all $i>1$. In particular, if $\delta(f)=0$, the integral $D$ can be taken with $\varphi_D(f)=f$.
\end{prop}

\section{Integrable derivations for $x^n-y^q$}

Let $R=k[x,y]$ be the polynomial ring in two variables over a reduced ring $k$ of characteristic $p>0$ and $h=x^n-y^q\in R$. In this section we will study the modules of $n$-integrable derivations of $A=R/\h$ of length $n\in \overline{\mathbb N}$.

In this section we will follow the following notation: Let $\alpha:=\val_p(n)$ be the p-adic valuation of $n$ and $s=n/p^\alpha$. We will denote by $m$ the remainder of the division of $q$ by $p$ and $\beta:=\val_p(q-m)$. We write $$\gamma:=\min\{i|ip^\alpha \geq q-1\}=\lceil (q-1)/p^\alpha\rceil.$$

\begin{prop}\label{integrabilidad binomio}
Let $k$ be a commutative reduced ring of characteristic $p>0$ and
$R=k[x,y]$ the polynomial ring over $k$. We set $A=R/\h$ where
$h=x^n-y^q$. For $\delta\in \Der_k(\log h)$, we denote
$\overline{\delta}=\Pi(\delta)$ (Corollary \ref{apli}).
\begin{itemize}
\item If $n,q\neq 0$, then
$$
\IDer_k(A)=\Der_k(A)=\langle
\overline{\delta_1},\overline{\delta_2}\rangle
$$
where $\delta_1=qx\partial_x+ny\partial_y$ and
$\delta_2=qy^{q-1}\partial_x+nx^{n-1}\partial_y$.

\item If $n=0\mod p$ and $q=1$, then
$$
\IDer_k(A)=\Der_k(A)= \langle \overline{\partial_x}\rangle
$$
\item If $\alpha,m\geq 1$ and $q\geq 2$, then
$$
\IDer_k(A;i)=\left\{ \begin{array}{ll} \left\{
\begin{array}{ll}
\left\langle \overline{\partial_x}\right\rangle& 1\leq i < p^\alpha\\
\left\langle\overline{ x\partial_x},\overline{ y^{\gamma}\partial_x}\right\rangle & p^\alpha\leq i< p^{\alpha+\beta}\\
\left\langle \overline{x\partial_x},
\overline{y^{\gamma+1}\partial_x}\right\rangle& i\geq
p^{\alpha+\beta} \mbox{ or } i=\infty
\end{array}
\right. & \mbox{ if } s=1,\mbox{ }\alpha\leq \beta,\mbox{ } m=1\\
\left\{
\begin{array}{lllll}
\langle \overline{\partial_x}\rangle &&\mbox{ }1\leq i<p^\alpha\\
\left\langle
\overline{x\partial_x},\overline{y^{\gamma}\partial_x}\right\rangle
&&\mbox{ } i\geq p^\alpha \mbox{ or } i=\infty
\end{array}\right. & \mbox{ otherwise}
\end{array}
\right.
$$
\end{itemize}

\end{prop}

\noindent{\bf Proof.}

Let $\delta=u\partial_x+v\partial_y$ be a $k$-derivation of $R$. To
prove this result it is enough to show which derivations are
$h$-logarithmically $i$-integrable for $i\in \overline{\mathbb{N}}$
(Corollary \ref{apli}).

\begin{itemize}
\item {\it $n,q\neq 0 \mod p$.}
\end{itemize}
We have to find the pairs $(u,v)$ such that $\delta(h)= nux^{n-1}-qvy^{q-1}\in \h$. It easy to see that $\Der_k(\log h)=\langle \delta_1,\delta_2\rangle$
where $\delta_1=qx\partial_x+ny\partial_y$ and
$\delta_2=qy^{q-1}\partial_x+nx^{n-1}\partial_y$. Note that $h$ is a quasi-homogenous polynomial with respect to the weights $w(x)=q$ and $w(y)=n$. By Theorem 1.2. of \cite{Tr}, the Euler vector field, $\delta_1$, is $h$-logarithmically $\infty$-integrable. On the other hand, the
gradient of $h$ is $J^0=\langle x^{n-1},y^{q-1}\rangle$, so
$\delta_2\in J^0\Der_k(R)$ and from Proposition \ref{2.2.4Na2}
we know that $\delta_2$ is $h$-logarithmically $\infty$-integrable
too. So,
$\IDer_k(A)=\Der_k(A)=\langle\overline{\delta_1}, \overline{\delta_2}\rangle.$

\begin{itemize}
\item {\it $n=0\mod p$ and $q=1$}.
\end{itemize}
The condition for $\delta$ to be $h$-logarithmic is that $v\in \h$,
so $\Der_k(\log h)=\langle \partial_x,h\partial_y\rangle$. In this
case $J^0=\langle 1\rangle$, hence any $\h$-logarithmic derivation
is integrable (Prop. \ref{2.2.4Na2}). Then,
$\IDer_k(A)=\Der_k(A)=\langle \overline\partial_x \rangle$.
\begin{itemize}
\item {\it $\alpha,m\geq 1$ and $q\geq 2$}.
\end{itemize}
Note that $n=sp^\alpha$. In order
for $\delta$ to be $h$-logarithmic, $v\in \h$ so
$\Der_k(\log h)=\left\langle \partial_x,h\partial_y\right\rangle$.
Since $h\partial_y$ is the zero derivation on $A$, we can focus on the $h$-logarithmically integrability of $\delta=u\partial_x$ with $u\in R$. Let $u_x\in k[x,y]$ and $u_y\in k[y]$ such that
$$
u=u_x(x,y)x+u_y(y) \Rightarrow
\delta=u\partial_x=u_xx\partial_x+u_y\partial_x.
$$
Since $h$ is a quasi-homogeneous polynomial with respect to the weights $w(x)=q$ and $w(y)=sp^\alpha$, the Euler vector field, $\chi=qx\partial_x$, is $h$-logarithmically integrable, and hence also $u_xx\partial_x$ are. Since $\IDer_k(\log h;i)$ is a $R$-modules for all $i$,
$$
\delta\in \IDer_k(\log h;i)\Leftrightarrow u_y\partial_x\in
\IDer_k(\log h;i)
$$
Let us consider $\delta=u\partial_x$ where $u\in k[y]$. Let $\varphi:R\rightarrow R[|\mu|]$ be a $k$-algebra homomorphism:
$$
\begin{array}{rccl}
\varphi:&R&\longrightarrow&R[|\mu|]\\
&x&\longmapsto& x+u\mu+u_2\mu^2+\cdots\\
&y&\longmapsto& y\mbox{ \hspace{2mm} }+ \mbox{ \hspace{2mm}
}v_2\mu^2+\cdots
\end{array}
$$
To show that $\delta$ is $i$-integrable it is enough to prove that there exist $u_j,v_j$ for $2\leq j\leq i$ such that $\varphi(h)\in \h \mod \mu^{i+1}$, or, equivalently, the coefficients of $\mu^j$ in $\varphi(h)$ belong to $\h$ for all $j\leq i$. We will denote by $\mu_j$ the coefficient of $\mu^j$ in the equation
\begin{equation}\label{laecu1}
\varphi(h)=\left(x^{p^{\alpha}}+u^{p^\alpha}\mu^{p^\alpha}+
u_2^{p^\alpha}\mu^{2p^\alpha}+\cdots\right)^{s}-
\left(y+v_2\mu^2+v_3\mu^3+\cdots\right)^q
\end{equation}

Suppose that there exists $i$ such that $2\leq i <p^\alpha$. Then,
$\mu_2=-q y^{q-1}v_2$ has to belong to $\h$. Hence, $v_2\in \h$, so
we can put $v_2=0$. Let us assume that $v_l=0$ for all
$2\leq l<i<p^{\alpha}$. In this case, $\mu_i=-qy^{q-1}v_i$ and, as the
same before, we can put $v_i=0$. Then,
$$
\Der_k(A)=\IDer_k(A;i)=\left\langle
\overline{\partial_x}\right\rangle \mbox{ } \forall i<p^\alpha
$$
and we can write the equation (\ref{laecu1}) as:
\begin{equation}\label{laecu2}
\left(x^{p^{\alpha}}+u^{p^\alpha}\mu^{p^\alpha}+
u_2^{p^\alpha}\mu^{2p^\alpha}+\cdots\right)^{s}-
\left(y+v_{p^\alpha}\mu^{p^\alpha}+
v_{p^{\alpha}+1}\mu^{p^\alpha+1}+\cdots\right)^q\in \h
\end{equation}
Now, we have to see that there are $u_{p^\alpha},v_{p^\alpha}\in R$ such that
\begin{equation}\label{v_1}
\mu_{p^\alpha}=
sx^{p^\alpha(s-1)}u^{p^\alpha}-qy^{q-1}v_{p^\alpha}\in \h
\end{equation}
 Since $u\in k[y]$, the previous expression implies that
$u^{p^\alpha}\in \langle y^{q-1}\rangle$. Therefore, if we write
$u=\sum_{i\geq 0}u_iy^i$ with $u_i\in k$, then $u_i^{p^\alpha}=0$ for all $i$ such that $ip^\alpha <q-1$, so
$u_i=0$ because $k$ is reduced. Hence, we can write $u=w(y)y^\gamma$ where $\gamma=\min\{i|ip^{\alpha} \geq q-1\}$
and $w(y)\in k[y]$. Substituting the expression of $u$ on
(\ref{v_1}), we can deduce that
\begin{equation}\label{Vpa}
sx^{p^\alpha(s-1)}w^{p^\alpha}y^{\gamma
p^{\alpha}-(q-1)}-qv_{p^\alpha}\in \h \mbox{  }\Rightarrow \mbox{  } v_{p^\alpha}\in
(s/q)x^{p^\alpha(s-1)}w^{p^\alpha}y^{\gamma p^{\alpha}-(q-1)}+\h
\end{equation}
Therefore, $A$ has a leap on $p^\alpha$ and
$$
\IDer_k(A; p^{\alpha})=\langle
\overline{x\partial_x},\overline{y^\gamma
\partial_x}\rangle \mbox{ where } \gamma=\min\{i|\mbox{
}ip^{\alpha}\geq q-1\}.
$$

Let us write $q=tp^\beta+m$.  Note that the only case where $\gamma p^\alpha = q-1$ is $q=tp^\beta +1$ and $\alpha \leq \beta$. Let us focus on this case when $s=1$.

\begin{itemize}
\item {\it  Case $q=tp^\beta +1$, $\alpha \leq \beta$ and $s=1$}. Observe that $t\neq 0$ because $q\geq 2$.   It is easy to see that $\gamma=tp^{\alpha-\beta}$. We will study the integrability of
 $w(y)y^\gamma\partial_x$ in this particular case.

Substituting the values of $q$ and $s$ in the equation
(\ref{laecu2}) and (\ref{Vpa}) we obtain:

$$
\begin{array}{l}
\left(x^{p^{\alpha}}+u^{p^\alpha}\mu^{p^\alpha}+
u_2^{p^\alpha}\mu^{2p^\alpha}+\cdots\right)-
\left(y^{p^\beta}+v_{p^{\alpha}}^{p^\beta}\mu^{p^{\alpha+\beta}}+
v_{p^{\alpha}+1}^{p^\beta}\mu^{(p^\alpha+1)p^\beta}+
\cdots\right)^{t}
\left(y+v_{p^\alpha}\mu^{p^\alpha}+\cdots\right)\in \h
\end{array}
$$
and
$$
v_{p^\alpha}= cw^{p^\alpha}+Fh
$$
for $c=1/q$ and some $F\in k[x,y]$. Let us consider $i$ such that
$p^\alpha< i<p^{\alpha+\beta}$. If $i=jp^\alpha$ for some $j\geq 2$, then
$\mu_i=u_j^{p^\alpha}-y^{tp^\beta}v_i$. Otherwise,
$\mu_i=-y^{tp^\beta}v_i$. So, $w y^\gamma\partial_x$ is
$h$-logarithmically $i$-integrable for all $i<p^{\alpha+\beta}$
(it's enough to put $u_j=v_i=0$ so that $\mu_i\in \h$). Now,
$$
\mu_{p^{\alpha+\beta}}=u_{p^\beta}^{p^\alpha}-ty^{(t-1)p^\beta+1}v_{p^\alpha}^{p^\beta}-
y^{tp^\beta}v_{p^{\alpha+\beta}}
$$
has to belong to $\h$. So, substituting the value of $v_{p^\alpha}$, we have that
$$
u_{p^\beta}^{p^\alpha}-ctw^{p^{\alpha+\beta}}y^{(t-1)p^\beta+1}-
y^{tp^\beta}v_{p^{\alpha+\beta}}=
G\left(x^{p^\alpha}-y^{tp^\beta+1}\right)
$$
for some $G\in k[x,y]$. The coefficient of $y^j$ with $j=(t-1)p^\beta+1$ in this equality is
$tcw_0^{p^\alpha}=0$ where $w_0$ is the independent term of $w$.
Since $R$ is reduced, $w_0=0$. Hence, $y^\gamma\partial_x$ is not
$p^{\alpha+\beta}$-integrable. However, if $w=w'y$ with $w'\in
k[y]$, the previous equation is
$$
\begin{array}{cc}
u_{p^\beta}^{p^\alpha}-ctw'^{p^{\alpha+\beta}}
y^{q+p^\beta(p^{\alpha}-1)}- y^{tp^\beta}v_{p^{\alpha+\beta}}=
G\left(x^{p^\alpha}-y^{tp^\beta+1}\right)
\end{array}
$$
Then, there exists a solution, for instance  $u_{p^\beta}=0$ and
$v_{p^{\alpha+\beta}}=-ctw'^{p^{\alpha+\beta}}y^{p^\beta(p^\alpha-1)+1}$.
In conclusion, in this case $A$ has a leap in $p^{\alpha+\beta}$ and
$$
\IDer_k\left(A;p^{\alpha+\beta}\right)=\left\langle
\overline{x\partial_x},\overline{y^{\gamma+1}\partial_x}\right\rangle
$$
\end{itemize}
Until now we saw that, for all $q\geq 2$
$$
\IDer_k\left(A;p^\alpha\right)=\left\langle
\overline{x\partial_x},\overline{y^\gamma\partial_x}\right\rangle \mbox{ where
} \gamma=\min\{i|\mbox{ }ip^\alpha\geq q-1\}
$$
and moreover, when $q=tp^\beta+1$, $1\leq \alpha
\leq \beta$ and $s=1$, $y^\gamma\partial_x$ is not
$h$-logarithmically integrable but
$$
\IDer_k\left(A;p^{\alpha+\beta}\right)=\left\langle
\overline{x\partial_x},\overline{y^{\gamma+1}\partial_x}\right\rangle
$$
Let us rewrite $\gamma:=\gamma+1$ in the latter case. We will see
that $y^\gamma\partial_x$ is integrable on $A$ for all $q\geq 2$. Consider
$$
\begin{array}{rccl}
\varphi:&A&\longrightarrow& A[|\mu|]\\
&x&\longmapsto&x+y^\gamma\mu\\
&y&\longmapsto&y+v_1\mu^{p^\alpha}+v_2\mu^{2p^\alpha}+\cdots
\end{array}
$$
where
$$
v_i=C_i x^{p^\alpha(s-\sigma)}y^{i\gamma p^{\alpha}-(\tau+1)q+1}
\mbox{ for } i=\tau s+\sigma \mbox{ with } \tau\geq 0 \mbox{ and }
\sigma=1,\ldots,s,
$$
$$
C_i=\frac{1}{q}\left[\binom{s}{i}-\sum_{j\in I_i}D_j\right] \mbox{
where } \binom{s}{i}=0 \mbox{ if } i>s,
$$
$$
I_i=\left\{j=\left(j_0,j_1,\ldots,j_{i-1}\right)|\mbox{ } j_k\geq 0
\mbox{ } \forall k=0,\ldots,i-1, |j|=q,
\sum_{k=1}^{i-1}kj_k=i\right\}
$$
and, for all $j=(j_0,j_1,\ldots,j_l)$ with $l\geq 1$,
$$
D_j=\binom{q}{j}C_1^{j_1}\cdots C_{l}^{j_{l}} \mbox{ with
}\binom{q}{j}=\dfrac{q!}{j_0!\cdots j_{l}!}.
$$
We have to prove that $\varphi$ is well defined. First we see that $i\gamma
p^\alpha-(\tau+1)q+1\geq 0$, i.e., $(\tau s +\sigma)\gamma p^\alpha -\tau q \geq q-1$.

\begin{itemize}
\item When $\gamma p^\alpha >q-1$, then $\gamma p^\alpha \geq q$, but
$q$ is not multiple of $p$, so $\gamma p^\alpha \geq q+1$ and
therefore
$$
(\tau s+\sigma)\gamma p^\alpha -\tau q \geq (\tau
s+\sigma)(q+1)-\tau q=(\tau(s-1)+\sigma)q+\tau s+\sigma\geq q-1
$$
because $s-1\geq 0$ and $\sigma \geq 1$.

\item Let us consider $\gamma p^\alpha=q-1$. As we have seen before, the previous equality only hold if $q=tp^\beta+1$ and $\alpha\leq \beta$. If $s=1$, then we have considered $\gamma+1$, so we are in the first point. Therefore, we have just considered $s\geq 2$. In this case, we have to prove that
$ (\tau s +\sigma)\gamma p^\alpha -\tau q=(\tau s +\sigma)(q-1)
-\tau q\geq q-1 $. Then
$$
(\tau s +\sigma)(q-1) -\tau q \geq (2\tau  +\sigma)(q-1) -\tau
q=(\tau +\sigma)q-(2\tau+\sigma)
$$
So,
$$
(\tau +\sigma)q-(2\tau+\sigma)  \geq q-1 \Leftrightarrow
(\tau+\sigma-1)q\geq 2\tau+\sigma -1
$$
and this is true because $q\geq 2$ and $\tau+\sigma -1 \geq 0$. Note
that if $\tau+\sigma -1=0$ then $\tau=0$ and $\sigma=1$, so
$2\tau+\sigma-1=0$ too.
\end{itemize}

Now, we have to show that $\varphi(h)=0$ in $A[|\mu|]$. The equation is:
$$
\varphi(h)=\left(x^{p^\alpha}+y^{\gamma p^\alpha}\mu^{p^\alpha}\right)^s-
\left(y+v_1\mu^{p^\alpha}+v_2\mu^{2p^\alpha}+\cdots\right)^q
$$
Since all degrees of the monomial which appeared in this equation
are multiple of $p^\alpha$, let us denote $\mu_i$ to the coefficient
of degree $ip^\alpha$. Then
$$
\mu_i= \binom{s}{i} x^{p^\alpha(s-i)} y^{i\gamma p^\alpha}
-\widetilde{\mu_i}
$$
where $\widetilde{\mu_i}$ is the coefficient of $\mu^{ip^\alpha}$
from $\left(y+v_1\mu^{p^\alpha}+v_2\mu^{2p^\alpha}+\cdots\right)^q$.
This coefficient can be found on
$$
\left(y+v_1\mu^{p^\alpha}+\cdots+v_i\mu^{ip^\alpha}\right)^q=
\sum_{|j|=q} \binom{q}{j}y^{j_0}v_1^{j_1}\cdots
v_{i}^{j_i}\mu^{p^\alpha(j_1+\ldots+ij_i)}
$$
We just have to consider all $j$ such that $j_1+\ldots+ij_i=i$.
Observe that there exists only one $j$ holding
this equation such that $j_i\neq 0$, This $j$ is $(q-1,0,\ldots,0,1)$ where 1 is in the position
$i$. So, we can identify the set of all these $j$ with $I_i\cup (q-1,0,\ldots,0,1)$. Let us calculate a term of $\widetilde{\mu_i}$. Fixed $j$, we have
$$
\binom{q}{j}y^{j_0}v_1^{j_1}\cdots v_i^{j_i}=\binom{q}{j}
C_1^{j_1}\ldots C_i^{j_i} x^{ap^\alpha} y^b=D_j x^{ap^\alpha}y^b
$$
where
$$
\begin{array}{lcc}
\displaystyle a=\sum_{1\leq \tau s+\sigma \leq i} j_{\tau s+
\sigma}(s-\sigma)\geq 0 & \mbox{ and }& \displaystyle
b=j_0+\sum_{1\leq \tau s+\sigma \leq i} j_{\tau s+
\sigma}\left(\gamma p^\alpha(\tau s+\sigma)-(\tau+1)q+1\right)\geq 0
\end{array}
$$
We are going to study these exponents.
$$
a=s\sum_{1\leq \tau s+\sigma \leq i} j_{\tau s+ \sigma}-\sum_{1\leq
\tau s+\sigma \leq i} j_{\tau s+ \sigma}\sigma=s(q-j_0)-\sum_{1\leq
\tau s+\sigma \leq i} j_{\tau s+ \sigma}\sigma
$$
On the other side, we have
$$
\sum_{1\leq \tau s+\sigma \leq i} j_{\tau s+ \sigma}(\tau
s+\sigma)=s\sum_{1\leq \tau s+\sigma \leq i} j_{\tau s+
\sigma}\tau+\sum_{1\leq \tau s+\sigma \leq i} j_{\tau s+
\sigma}\sigma=ls+r
$$
where $i=ls+r$ (remember: $l\geq 0$ and $1\leq r\leq s$). Then, if we denote $T= \sum\limits_{1\leq \tau s+\sigma \leq i} j_{\tau s+
\sigma}\tau$ and we substitute on $a$, we have
$$
a= s(q-j_0)-((l-T)s+r)=s(q-j_0-l+T)-r\geq 0
$$
If $q-j_0-l+T< 1$, then $a<0$ so $q-j_0-l+T\geq 1$ and we can write
$$
a=(q-j_0-l+T-1)s+s-r
$$
Observe that  $s-r\geq 0$ because $1\leq r\leq s $. Now,
$$
\begin{array}{rll}
 b=&\displaystyle j_0+\gamma p^\alpha\sum_{1\leq \tau s+\sigma \leq i} j_{\tau s+ \sigma}(\tau s+\sigma)-q\sum_{1\leq \tau s+\sigma \leq i} j_{\tau s+ \sigma}\tau-q\sum_{1\leq \tau s+\sigma \leq i} j_{\tau s+ \sigma}+\sum_{1\leq \tau s+\sigma \leq i} j_{\tau s+ \sigma}\\[8mm]
=&\gamma p^\alpha i-qT-q(q-j_0)+(q-j_0)+j_0=i\gamma p^\alpha
-q(T+q-j_0-1)
\end{array}
$$
So,
$$
\binom{q}{j}y^{j_0}v_1^{j_1}\cdots
v_i^{j_i}=D_jx^{(q-j_0-l+T-1)sp^\alpha+(s-r)p^\alpha} y^{i\gamma
p^\alpha-q(T+q-j_0-1)}
$$
Since $x^{sp^\alpha}=y^q$ in $A$,
$$
\binom{q}{j}y^{j_0}v_1^{j_1}\cdots v_i^{j_i}=D_j x^{(s-r)p^\alpha}
y^{i\gamma
p^\alpha+q(q-j_0-l+T-1)-q(T+q-j_0-1)}=D_jx^{(s-r)p^\alpha}y^{i\gamma
p^\alpha-lq}
$$
Hence,
$$
\begin{array}{rll}
\widetilde{\mu_i}=&\displaystyle \sum_{\substack{|j|=q\\
j_1+\ldots+ij_i=i}} D_jx^{p^\alpha(s-r)}y^{i\gamma p^\alpha-lq}=
\left(\sum_{j\in I_i} D_j+D_{(q-1,0,\ldots,0,1)}\right)x^{p^\alpha(s-r)}y^{i\gamma p^\alpha-lq}\\
=&\displaystyle \left(\sum_{j\in I_i}
D_j+qC_i\right)x^{p^\alpha(s-r)}y^{i\gamma
p^\alpha-lq}=\left(\sum_{j\in I_i}
D_j+q(1/q)\left[\binom{s}{i}-\sum_{j\in
I_i} D_j\right] \right) x^{p^\alpha(s-r)}y^{i\gamma p^\alpha-lq}\\
\displaystyle =&\displaystyle \binom{s}{i}
x^{p^\alpha(s-r)}y^{i\gamma p^\alpha-lq}
\end{array}
$$
So,
$$
\mu_i=\binom{s}{i} x^{p^\alpha(s-i)}y^{i\gamma p^\alpha}-\binom{s}{i}x^{p^\alpha(s-r)}y^{i\gamma p^\alpha-lq}
$$

If $i>s$, then $\binom{s}{i}= 0$, and hence $\mu_i=0$.

If $i\leq s$, then $i=0 \cdot s+i$, i.e., $l=0$ and $r=i$, then
$$
\mu_i=\binom{s}{i} x^{p^\alpha(s-i)} y^{i\gamma
p^\alpha} -\binom{s}{i} x^{p^\alpha(s-i)}y^{i\gamma p^\alpha}=0
$$
so, $\varphi$ is well defined and the proposition is proved.
\begin{flushright}$\square$\end{flushright}

\begin{ej}
Let us consider $k$ a reduced ring of characteristic $p=3$ and $h=x^3-y^4\in k[x,y]$,
then $\gamma=1$ so, according with Proposition \ref{integrabilidad binomio},
$$
\IDer_k(A;i)=\left\{\begin{array}{ll} \langle \overline\partial_x\rangle& 1\leq i<3\\
\langle \overline{x\partial_x}, \overline{y\partial_x}\rangle& 3\leq i<9\\
\langle \overline{x\partial_x},\overline{y^2\partial_x}\rangle& i\geq 9
\end{array}
\right.
$$
Now, if we consider $h=x^3-y^5$, then $\gamma=2$ and
$$
\IDer_k(A;i)=\left\{\begin{array}{ll} \langle \overline\partial_x\rangle& 1\leq i<3\\
\langle \overline{x\partial_x}, \overline{y^2\partial_x}\rangle& i\geq 3\\
\end{array}
\right.
$$
\end{ej}

\begin{nota}
Note that if $k$ is not reduced, Proposition \ref{integrabilidad binomio} is not true. For example, if $k=\mathbb F_3[t]/\langle t^3\rangle$ and $h=x^3-y^5$, then $\overline{t\partial_x}\in \IDer_k(A)$ with the integral
$$
\begin{array}{ccccccccccccccccccccccccccccccccccccccc}
A&\rightarrow&A[|\mu|]\\
x&\mapsto& x+t\mu\\
y&\mapsto& y
\end{array}
$$
\end{nota}

\begin{cor}
Let $k$ be a commutative reduced ring of characteristic $p>0$ and $A=k[x,y]/\h$ where $h=x^n-y^q$ such that $\alpha,m\geq 1$ and $q\geq 2$. We denote  $B_i:=\Ann_A\left(\IDer_k(A;i-1)/\IDer_k(A;i)\right)$ for $i>1$.
Then,
$$
B_i=\left\{ \begin{array}{lll}
\langle x,y^\gamma\rangle &\mbox{ if } i=p^\alpha\\
\langle y\rangle & \mbox{ if } i=p^{\alpha+\beta}, \mbox{ } s=1,\mbox{ } \alpha\leq \beta \mbox{ and } m=1\\
\end{array}
\right.
$$
Moreover, $B_i\supseteq J^0=\langle y^{q-1}\rangle$ where $J^0$ is the gradient ideal of $h$ defined in Proposition \ref{2.2.4Na2}.
\end{cor}

\noindent{\bf Proof.}

Let us start with $i=p^\alpha$. From Proposition \ref{integrabilidad
binomio}, we can deduce that
$$
\IDer_k\left(A;p^\alpha-1\right)/\IDer_k\left(A;p^\alpha\right)=
\langle \partial_x\rangle/\langle x\partial_x,y^\gamma\partial_x\rangle
$$
where $\partial_x\in \Der_k(A)$.  By definition, $a\in B_i$ if
$a\partial_x=0 \mod \langle x\partial_x,y^\gamma\partial_x\rangle$,
i.e, if there exist $F,G\in A$ such that
$a\partial_x=Fx\partial_x+Gy^\gamma
\partial_x$. Applying this derivation to $x$, we have that $a\in
\langle x,y^\gamma \rangle$.

Now, when $\alpha\leq \beta$,  $s=m=1$ and $i=p^{\alpha+\beta}$,
from Proposition \ref{integrabilidad binomio},
$$
\IDer_k\left(A;p^{\alpha+\beta}-1\right)/\IDer_k\left(A;p^{\alpha+\beta}\right) =\langle x\partial_x,y^\gamma\partial_x\rangle /\langle x\partial_x,y^{\gamma+1}\partial_x\rangle=\langle y^\gamma\partial_x/y^{\gamma+1}\partial_x\rangle
$$
In this case, $a\in B_{p^\alpha+\beta}$ if and only if
$ay^\gamma\partial_x\in\langle y^{\gamma+1}\partial_x \rangle$, i.e,
if $(a-Fy)y^\gamma\partial_x=0$ for some $F\in A$.  This implies that $a\in \langle y\rangle$ and we have proved the corollary.
\begin{flushright}$\square$\end{flushright}

\subsection{$I^p$-logarithmic derivations}

In this section, we want to calculate the $m$-integrable derivations of
$A=k[x,y]/\h$ where $k$ is a unique factorization domain (UFD) of
characteristic $p>0$ and $h=x^n-y^q$ with $n,q=0\mod p$. We start
with some general results about the relationship between $\langle
f\rangle$-logarithmic and $\langle f^p\rangle $-logarithmic
derivations. In this section, we denote $R=k[x_1,\ldots,x_d]$.

\begin{prop}\label{HS de producto primo}
If $f,g\in R=k[x_1,\ldots,x_d]$ are coprime, then, for all $n\in \overline{\mathbb
N}$, we have:
$$
\HS_k(\log fg;n)=\HS_k(\log f;n)\cap \HS_k(\log g;n).
$$
\end{prop}

\noindent{\bf Proof.}
\begin{itemize}
\item[$\supseteq$.] Let $D\in \HS_k(\log f;n)\cap \HS_k(\log g;n)$. By definition, $D_i(f)\in \langle f\rangle$ and $D_i(g)\in \langle
g\rangle$ for all $i\leq n$. Then $D_i(fg)=\sum_{a+b=i}
D_a(f)D_b(g)\in \langle fg\rangle$, so $D\in \HS_k(\log fg;n)$.

\item[$\subseteq$.] Let $D\in \HS_k(\log fg;n)$. This implies that $D_i(fg)\in
\langle fg\rangle$ for all $i\leq n$. We will prove the result by
induction on $i$. When $i=1$, then $D_1(fg)=D_1(f)g+D_1(g)f\in
\langle fg \rangle \subseteq \langle f\rangle, \langle g \rangle$.
So, $D_1(f)g\in \langle f \rangle$. Since $g$ and $f$ are coprime,
$D_1(f)\in \langle f \rangle$. For $g$ is analogous.

Now let us assume that $D_i(f)\in \langle f\rangle$ and $D_i(g)\in
\langle g \rangle$ for all $i<n$. By definition,
$$
D_n(fg)=D_n(f)g+D_n(g)f +\sum_{\substack {a+b=n\\ a,b\neq 0}}
D_a(f)D_b(g)\in \langle fg\rangle \Rightarrow D_n(f)g+D_n(g)f\in
\langle fg\rangle
$$
and we can proceed like case $i=1$.
\end{itemize}
\vspace{-1.0cm}\begin{flushright}$\square$\end{flushright}

\begin{cor}\label{IDer de producto primo}
If $f,g\in R$ are coprime, then $ \IDer_k(\log fg;n)\subseteq
\IDer_k(\log f; n)\cap \IDer_k(\log g;n)$ for all $n\in
\overline{\mathbb N}$.
\end{cor}

\noindent{\bf Proof.} If $\delta\in \IDer_k(\log fg;n)$ then, there
exists $D\in \HS_k(\log fg;n)$ a $n$-integral of $\delta$. By Proposition \ref{HS de producto primo}, $D\in
\HS_k(\log f;n)\cap\HS_k(\log g;n)$ so, $\delta\in
\IDer_k(\log f;n)\cap \IDer_k(\log g;n)$.
\begin{flushright}$\square$\end{flushright}

\begin{nota}
In general, equality in Proposition 2.5 does not hold. For example: Let $k=\mathbb F_2$
and $f=y^2$ and $g=x^2-y$ two polynomial of $k[x,y]$. Then $
\partial_x\in \IDer_k(\log f;4)\cap \IDer_k(\log g;4)
$. However $\partial_x \not\in \IDer_k(\log fg;4)$.
\end{nota}

\begin{cor}\label{Producto de varios primos}
Let $f_1,\ldots, f_m\in R$. If $f_i$,$f_j$ are coprime whenever
$i\neq j$, then, for all $\overline{\mathbb N}$ we have:
$$
\begin{array}{ccc}
\HS_k(\log f_1\cdots f_m;n)=\bigcap_i \HS_k(\log f_i;n) &\mbox{and
}& \IDer_k(\log f_1\cdots f_m; n)\subseteq
\bigcap_i\IDer_k(\log f_i;n)
\end{array}
$$
\end{cor}

\noindent{\bf Proof.} The result is obtained thanks to Proposition
\ref{HS de producto primo} and Corollary \ref{IDer de producto
primo} by induction on $m$.
\begin{flushright}$\square$\end{flushright}

\begin{lem}\label{Rel1}
Let $f$ be an irreducible polynomial, $a\geq 1$ and $n\in
\overline{\mathbb N}$. Consider $D\in \HS_k(R;n)$. Suppose that
$D_i(f^a)^p\in \left\langle f^{ap}\right\rangle$ for all $i\leq n$.
Then, $D\in \HS_k(\log f^a;n)$.
\end{lem}

\noindent{\bf Proof.}

We write $a=sp^\alpha$ where $\alpha=\val_p(a)\geq 0$ and $s\geq 1$. By Lemma \ref{Dihp},
$$
D_i\left(f^{s p^\alpha}\right)=\left\{\begin{array}{lcl}
0&\mbox{ if }& p^\alpha \not| i\\
D_{i/p^\alpha}(f^s)^{p^\alpha}& \mbox{ if }&p^\alpha|i
\end{array}\right.
$$
Hence, we can focus on the case $n\geq p^\alpha$ and $i=jp^\alpha\leq n$. It's enough to show that $D_j(f)\in \langle f\rangle$ because, if this is true, we have that $D_j(f^s)\in \langle f^s\rangle$ by Lemma \ref{ExpDi},
and $D_i\left(f^{sp^\alpha}\right)=D_j(f^s)^{p^\alpha}\in
\left\langle f^{sp^\alpha}\right\rangle$ so we would have the result.

\vspace{0.2cm}
Since $i=jp^\alpha\leq n$,
\begin{equation}\label{Hip}
D_j\left(f^s\right)^{p^{\alpha+1}}=D_{jp^\alpha}
\left(f^{sp^\alpha}\right)^p\in \left\langle
f^{sp^{\alpha+1}}\right\rangle
\end{equation}
When $j=1$, $D_1\left(f^s\right)=s f^{s-1}D_1(f)$. Substituting in the previous expression, we have that
\begin{equation}\label{Pert1}
D_1\left(f^s\right)^{p^{\alpha+1}}=s
f^{(s-1)p^{\alpha+1}}D_1(f)^{p^{\alpha+1}}\in \left\langle f^{s
p^{\alpha+1}}\right\rangle
\end{equation}
Since $R$ is UFD and  $f,s\neq 0$, $D_1(f)^{p^{\alpha+1}}\in
\left\langle f^{p^{\alpha+1}}\right\rangle\subseteq \langle f
\rangle$ and hence $D_1(f)\in \langle f \rangle$.

Let us assume that $D_l(f)\in \langle f\rangle $ for all $l<j$ with
$jp^\alpha\leq n$. Thanks to the hypothesis, we can use Lemma
\ref{ExpDi}, and we have
    $$
    D_j\left(f^s\right)=s f^{s-1}D_j(f)+Ff^s
    $$
for some $F\in R$. Substituting this expression in  (\ref{Hip}),
$$
s f^{(s-1)p^{\alpha+1}}D_j(f)^{p^{\alpha+1}}+F^{p^{\alpha+1}} f^{s
p^{\alpha+1}}\in \left\langle f^{s p^{\alpha+1}}\right\rangle
\Rightarrow s f^{(s-1)p^{\alpha+1}}D_j(f)^{p^{\alpha+1}}\in
\left\langle f^{s p^{\alpha+1}}\right\rangle
$$
Observe that it is the same condition that (\ref{Pert1}), so we can
deduce that  $D_j(f)\in \langle f\rangle$.
\vspace{-0.5cm}\begin{flushright}$\square$\end{flushright}

\begin{prop}\label{h-hp}
Let $k$ be an UFD of characteristic $p>0$ and $R=k[x_1,\ldots, x_d]$ the polynomial
ring over $k$.  Let $h$ be a polynomial of $R$. For all $n\in
\overline{\mathbb N}$, we have:
$$
\IDer_k(\log h;n)=\IDer_k\left(\log h^p,np\right).
$$
\end{prop}

\noindent{\bf Proof.}

$\subseteq$. Let $D_1\in \IDer_k(\log h;n)$ and $D\in \HS_k(\log h;n)$
an integral. If $n<\infty$, from Corollary \ref{HSenpolise
extiende},  $D$ is $np$-integrable, so let $D'$ be a $np$-integral
of $D$. If $n=\infty$, we put $D'=D$. Observe that $D'_1=D_1$ so, if
$D'\in \HS_k\left(\log  h^p;np\right)$ then $D_1\in
\IDer_k(\log h^p;np)$. We have to see that $D'_i(h^p)\in \langle h^p
\rangle$ for all $i\leq np$.

By Lemma \ref{Dihp},
$$
D_i'\left(h^p\right)=\left\{\begin{array}{lcl}
0&\mbox{ if }& p \not| i\\
D'_{i/p}(h)^{p}& \mbox{ if }&p|i
\end{array}\right.
$$
Then, we can focus on $i=jp$ where $1\leq j\leq n$. Note that
$D_j'=D_j$ for all $1\leq j\leq n$, so
$$
D_i'\left(h^p\right)=D'_j(h)^p=D_j(h)^p\in \langle h^p\rangle.
$$
Therefore, $D'_i(h^p)\in \langle h^p\rangle$ for all $i\leq np$ and we have
the inclusion.

$\supseteq$. Let $D_1\in \IDer_k\left(\log h^p;np\right)$ and $D\in
\HS_k\left(\log h^p;np\right)$ a $np$-integral of $D_1$. Let
$h=h_1^{a_1}\cdots h_m^{a_m}$ be the factorization of $h$ in
irreducible factors, i.e, $h_i$ is irreducible and $a_i\geq 1$ for all
$i=1,\ldots,m$ and $h_i\neq h_j$ if $i\neq j$. Then $h_i^{a_i}$ and $h_j^{a_j}$
are coprime whenever $i\neq j$, and therefore, $h_1^{a_1p},\ldots,
h_m^{a_mp}$ are coprime too. By Corollary \ref{Producto de
varios primos},
$$
D\in \HS_k(\log  h^p;np)= \bigcap_i \HS_k\left(\log  h_i^{a_ip};np\right).
$$
Hence, $D_j(h_i^{a_i})^p=D_{jp}(h_i^{a_ip})\in \langle h_i^{a_ip}
\rangle $ for $j\leq n$. By Lemma \ref{Rel1}, $D_j(h_i^{a_i})\in
\langle h_i^{a_i}\rangle$ for all $i=1,\ldots,m$, and $j\leq n$. So,
$ \tau_{np,n}(D)\in \cap\HS_k(\log h_i^{a_i};n)=\HS_k(\log h;n)
\Rightarrow D_1\in \IDer_k(\log h;n)$.
\vspace{-0.5cm}\begin{flushright}$\square$\end{flushright}

\begin{cor}\label{IDerhhpt}
For all $\tau \geq 0$ and $n\in \overline{\mathbb N}$, $
\IDer_k(\log h;n)=\IDer_k\left(\log h^{p^\tau};np^\tau\right) $.
\end{cor}

\noindent{\bf Proof.} By induction on $\tau$ using Proposition
\ref{h-hp}.
\vspace{-0.8cm}\begin{flushright}$\square$\end{flushright}

\begin{prop}\label{Saltos de hpt}
Let $k$ be a UFD of characteristic $p>0$, $R=k[x_1,\ldots,x_d]$ the polynomial ring
over $k$, $h\in R$ and $\tau\geq 1$.  Then the set of the leaps of
$A:=R/\left\langle h^{p^\tau}\right\rangle$ is
$$
\left\{\begin{array}{lcc}
\{np^\tau |\mbox{ }n \mbox{ leap of } R/\h \} &\mbox{ if } \Der_k\left(\log h\right)=\Der_k(R)\\
\{np^\tau |\mbox{ }n \mbox{ leap of } R/\h \} \cup {p^\tau} &\mbox{ if } \Der_k\left(\log h\right)\neq\Der_k(R)\\
\end{array}
\right.
$$
\end{prop}

\noindent{\bf Proof.}

By Corollary \ref{Saltos del cociente}, $A$ has a leap on $s>1$
if and only if  the inclusion $
\IDer_k\left(\log  h^{p^\tau};s-1\right)\supsetneq
\IDer_k\bigl(\log  h^{p^\tau};s\bigr) $ is proper. First of all, we will prove the next two equalities:
\begin{itemize}
\item[1.] For $s<p^\tau$, $\IDer_k\left(\log  h^{p^\tau};s\right)=\Der_k(R)$.
\end{itemize}
$\subseteq$ is always true. Let $D_1\in \Der_k(R)=\IDer_k(R)$ and
$D\in \HS_k(R)$ an integral. Since $s<p^\tau$, for all $j\leq s$,
$p^\tau \nmid j$. By Lemma \ref{Dihp},
$D_j\left(h^{p^\tau}\right)=0\in \left\langle
h^{p^\tau}\right\rangle$ for all $j\leq s$. Then, any derivation
$D_1$ has a $h^{p^\tau}$-logarithmic $s$-integral and the other
inclusion holds. So, $A$ does not have a leap on $s$.
\begin{itemize}
\item[2.] Let $s$ be an integer such that $np^\tau <s<(n+1)p^\tau$ for some $n\geq 1$.
Then $ \IDer_k\left(\log h^{p^\tau}; s\right)=
\IDer_k\left(\log h^{p^\tau}; np^\tau\right) $.
\end{itemize}
Since $s>np^\tau$, the inclusion $\subseteq$ is true. Let
$D_1\in\IDer_k\left(\log h^{p^\tau}; np^\tau\right)$. By definition
there exists an integral $D\in
\HS_k\left(\log h^{p^\tau};np^\tau\right)$. By Corollary \ref{HSenpolise extiende}, we can consider $D'\in \HS_k(R;s)$ an integral of $D$. Hence, for all $j$ such that $np^\tau<j\leq s<(n+1)p^\tau$, $p^\tau \not|j$ and, by Lemma \ref{Dihp}, $D'_j\left(h^{p^\tau}\right)=0 \in \left\langle
h^{p^\tau}\right\rangle $. Since $D'_l=D_l$ for all $l\leq np^\tau$, $D'\in
\HS_k\left(\log h^{p^\tau};s\right)$. Therefore, $D_1\in
\IDer_k\left(\log h^{p^\tau}; s\right)$ and $A$ does not have a leap
on $s$.

\vspace{0.2cm}
Thanks to these two equalities we know that the leaps are given on
$s=np^\tau$ for some $n\geq 1$. Let us suppose that $s=p^\tau$. By
Corollary \ref{IDerhhpt} and the point 1.,
$$
\Der_k\left(R\right)=
\IDer_k\left(\log h^{p^\tau};s-1\right)\supseteq
\IDer_k\left(\log  h^{p^\tau};p^\tau=s\right)=\Der_k\left(\log h\right)
$$
Hence, $A$ has a leap on $p^\tau$ if and only if $\Der_k(\log h
)\neq
\Der_k(R)$. Now, let us consider $s=np^\tau$ for $n\geq 2$. By Corollary \ref{IDerhhpt} and the point 2.
    $$
    \begin{array}{rlll}
\IDer_k\left(\log h; n-1\right)=& \IDer_k\left(\log h^{p^\tau};
(n-1)p^\tau\right) =  \IDer_k\left(\log h^{p^\tau};
np^\tau-1\right)\\ \supseteq&\IDer_k\left(\log h^{p^\tau};
np^\tau\right)=\IDer_k\left(\log h; n\right)
   \end{array}$$
Then, $A$ has a leap on $s=np^\tau$ if and only if $n$ is a leap on
$R/\h$ and we have proved the result.
\vspace{-0.5cm}\begin{flushright}$\square$\end{flushright}


\begin{prop}
Let $k$ be a UFD of characteristic $p>0$ and $h=x^n-y^q\in k[x,y]$. Suppose $\alpha:=\val_p(n)\geq 1$ and $\beta:=\val_p(q)\geq 1$. We write $\tau=\min{\alpha,\beta}\geq 1$, $s=n/p^\tau$ and $t=q/p^\tau$. Then,
$$
\IDer_k(k[x,y]/\h;np)=\left\{\overline{\delta}|\mbox{ } \delta\in \IDer_k(\log \langle x^s-y^t \rangle,n)\right\}
$$
where the leaps occur in $\left\{np^\tau|\mbox{ }n \mbox{ is a leap of } k[x,y]/\langle
H\rangle\right\} \cup {p^\tau}$.
\end{prop}

\noindent{\bf Proof.} Using Corollary \ref{IDerhhpt} and Proposition \ref{integrabilidad
binomio}. \vspace{-0.5cm}
\begin{flushright}$\square$\end{flushright}

\section{Other examples}

We are going to calculate the integrable derivations of the quotient of a polynomial ring over some non-binomial equations. These examples have been taken from the article \cite{Gr}.
\vspace{0.2cm}

\noindent{\it Example 1.}
\vspace{0.2cm}

Let $k$ be a domain of characteristic $p>0$ and $h=x^p+tx^{p+1}\in R=k[x]$ with $t\in k$. Let $A=R/\h$. The module of $\Der_k(\log h)$ is generated by $(1+tx)\partial_x$. From Example (2.1.2) of \cite{Na2}, we have that $(1+tx)\partial_x$ is $h$-logarithmically $(p-1)$-integrable. So, let us consider $E\in \HS_k(\log h;p-1)$ an integral of $u(1+tx)\partial_x$ where $u\in R$. From Corollary \ref{HSenpolise extiende}, there exists $D\in \HS_k(R)$ an integral of $E$. In order for $D$ to be $h$-logarithmic,
$$
D_p(x^p+tx^{p+1})=D_1(x)^p+t(xD_1(x)^p+D_p(x)x^p)= u^p(1+tx)^{p+1}+tD_p(x)x^p\in \h
$$
So, $u\in \langle x\rangle$ and $\IDer_k(\log h;p)=\langle x(1+tx)\partial_x\rangle$. Observe that this generator is $\infty$-integrable, for example $x\in A\mapsto x+x(1+tx)\mu\in A[|\mu|]$ is an integral. In conclusion,
$$
\IDer_k(A;i)=\left\{\begin{array}{lll}
\langle\overline{(1+tx)\partial_x} \rangle& \mbox{ if } i\leq p-1\\
\langle \overline{x(1+tx)\partial_x}\rangle &\mbox{ if } i\geq p
\end{array}
\right.
$$

\noindent{\it Example 2.}
\vspace{0.2cm}

Let $k$ be a domain of characteristic $p=2$ and $h=x^4+y^6+y^7\in R=k[x,y]$. Let $A=R/\h$. In this case, the module of $h$-logarithmic derivations is generated by $\partial_x$ and $h\partial_y$. Since $h\partial_y$ is $h$-logarithmically $\infty$-integrable, we can focus on the $h$-logarithmically integrability of $u\partial_x$ where $u\in k[x,y]$. Let $\varphi:R\rightarrow R[|\mu|]$ a $k$-algebra homomorphism:
$$
\begin{array}{rccc}
\varphi:&R&\longrightarrow&R[|\mu|]\\
&x&\longmapsto& x+u\mu+u_2\mu^2+\cdots\\
&y&\longmapsto& y\mbox{ \hspace{2mm} }+ \mbox{ \hspace{2mm}
}v_2\mu^2+\cdots
\end{array}
$$
We want to see that there exist $u_i,v_i\in R$ for $i\geq 2$ such that $\varphi$ is $h$-logarithmic.
The coefficient of $\mu^i$ for $i=2,3$ in $\varphi(h)$ is $y^6v_i$.
In order for $\varphi$ to be $h$-logarithmic, $v_i\in \h$, so we can
put $v_i=0$. In fact, we can put $v_i=0$ for all $i$ such that
$4\not| i$. Thanks to this, we can write:
\begin{equation}\label{1}
\varphi(h)=(x+u\mu+u_2\mu^2+\ldots)^4+ (y+v_4\mu^4+v_8\mu^8\ldots)^6(1+y+v_4\mu^4+v_8\mu^8\ldots)
\end{equation}

The coefficient of $\mu^4$ in (\ref{1}) is $\mu_4:=u^4+y^6v_4$ and it has to
belong to $\h$. Hence, $u\in \langle x,y^2\rangle$ and $\IDer_k(\log
h;4)=\langle x\partial_x,y^2\partial_x, h\partial_y\rangle$. It's
easy to proof the following lemma through the calculation of a term in the equation (\ref{1}):
\begin{lem}\label{Lema 1}
Suppose that $u_j=0$ for all $j\geq 2$ and $v_{4n}\in \langle
y^2\rangle$ for all $n<i$, then there exists $v_{4i}\in \langle
y^2\rangle$ such that the coefficient of $\mu^{4i}$ in (\ref{1}) belongs to $\h$.
\end{lem}

Using this lemma repeatedly we deduce that $y^2\partial_x$ and
$xy\partial_x$ are $h$-logarithmically integrable since a possible
solution so that $\mu_4$ is $h$-logarithmic is $v_4=y^2$ and $v_4
=(1+y)y^4$ respectively. Therefore, we need to see the $h$-logarithmically integrability
of $ux\partial_x$ where $u\in k[x]$. In this case, $v_4\in
(1+y)u^4+\h$. Calculating the coefficient of $\mu^8$ in (\ref{1}),
we obtain $\mu_8:=u_2^4+y^6v_8+v_4^2(1+y)y^4$. In order for $\mu_8$ to be in $\h$, $u\in
\langle x\rangle$. Hence, $v_4\in \langle x^4,h\rangle$. We deduce that
$x^2\partial_x$ is $h$-logarithmically integrable by the following
lemma:

\begin{lem}\label{Lema 2}
Suppose that $u_j=0$ for all $j\geq 2$ and  $v_{4n}\in \langle
x^4\rangle$ for all $n<i$, then there exists $v_{4i}\in \langle
x^4\rangle$ such that the coefficient of $\mu^{4i}$ in (\ref{1}) belongs to $\h$.
\end{lem}
In conclusion,
$$
\IDer_k(A;i)=\left\{ \begin{array}{lll}
\langle\overline{\partial_x}\rangle&\mbox{ if } 1\leq i<4\\
\langle \overline{x\partial_x},\overline{y^2\partial_x}\rangle& \mbox{ if } 4\leq i<8\\
\langle \overline{x^2\partial_x},\overline{xy\partial_x}, \overline{y^2\partial_x}\rangle& \mbox{ if } i\geq 8
\end{array}
\right.
$$

\noindent{\it Example 3.}
\vspace{0.2cm}

Let $k$ be a domain of characteristic $p=3$ and $h=x^3+y^5+x^2y^2\in
R=k[x,y]$.  Let $A=R/\h$. The module of $h$-logarithmic derivation
is generated by $\delta_1:=x^2\partial_x+y^3\partial_y$ and
$\delta_2:=2y^2\partial_x+(x+y^2)\partial_y$. These two derivations
are $h$-logarithmically integrable. To verify this claim, let us consider $\varphi:R\rightarrow R[|\mu|]$ a homomorphism of $k$-algebras
$$
\begin{array}{rccc}
\varphi:&R&\longrightarrow&R[|\mu|]\\
&x&\longmapsto& x+u_1\mu+u_2\mu^2+\cdots\\
&y&\longmapsto& y+v_1\mu+v_2\mu^2+\cdots
\end{array}
$$
As in the previous example, we want to prove that there exist $u_i,v_i\in R$ for $i\geq 2$ such that $\varphi(h)\in \h$ where $u_1$ and $v_1$ are determined by $\delta_1$ or $\delta_2$. By calculating a generic term of $\varphi(h)$, we
can show the following lemmas:
\begin{lem}\label{lema 3}
Let $u_1=x^2$ and $v_1=y^3$. Suppose that $v_j=0$ for all $j\geq 2$ and $u_n\in \langle x^2\rangle$ for all $n<i$. Then, there exists $u_i\in \langle x^2\rangle$ such that the coefficient of $\mu^i$ in $\varphi(h)$ belongs to  $\h$.
\end{lem}

\begin{lem}\label{lema 4}
Let $u_1=2y^2$ and $v_1=x+y^2$. Suppose $u_n\in \langle xy,y^3\rangle$ and $v_n\in \langle y^2\rangle$ for all $2\leq n<i$. Then, there exist $u_i\in \langle xy,y^3\rangle$ and $v_i\in \langle y^2\rangle$ such that the coefficient of $\mu^i$ in $\varphi(h)$ belongs to $\h$.
\end{lem}

Using Lemma \ref{lema 3} for the integrability of $\delta_1$ and Lemma \ref{lema 4} for the integrability of $\delta_2$, we can deduce that
$$
\IDer_k(A)=\langle \overline{\delta_1}, \overline{\delta_2}\rangle.
$$

\end{document}